\documentclass{article}
\usepackage[utf8]{inputenc}
\usepackage{amsfonts}
\usepackage{amssymb}
\usepackage{hyperref}
\usepackage{indentfirst}
\usepackage{mathrsfs}
\usepackage{tikz}
\usepackage{mathtools}
\usepackage{graphicx}
\usepackage{xparse}
\usepackage[all]{xy}

\usepackage[top=1in,bottom=1in,left=1.25in,right=1.25in]{geometry}
\date{}
\allowdisplaybreaks

\begin{document}

\renewcommand{\baselinestretch}{1.2}
\renewcommand{\arraystretch}{1.0}

\title{\bf Ore Extension of Group-cograded Hopf Coquasigroups}
\author
{\textbf{ Lingli Zhu} \footnote{College of Science, Nanjing Agricultural University, Nanjing 210095, Jiangsu, China.  E-mail: 2021111005@stu.njau.edu.cn}, \,
\textbf{Bingbing Jin} \footnote{College of Science, Nanjing Agricultural University, Nanjing 210095, Jiangsu, China.  E-mail: 2022111009@stu.njau.edu.cn }, \,
\textbf{Huili Liu} \footnote{Xixi branch school of Jiangsu Xishan Senior High School, Wuxi 214151, Jiangsu, China. E-mail: 710980586@qq.com }, \,
\textbf{Tao Yang} \footnote{Corresponding author. College of Science, Nanjing Agricultural University, Nanjing 210095, Jiangsu, China.     E-mail: tao.yang@njau.edu.cn}
}
\maketitle

\begin{center}
\begin{minipage}{12.cm}
 {Abstract: The aim of this paper is the Ore extension of group-cograded Hopf coquasigroups. This paper first shows a categorical interpretation and some examples of group-cograded Hopf coquasigroups, and then gives a necessary and sufficient conditions for the Ore extensions of group-cograded Hopf coquasigroups to be group-cograded Hopf coquasigroups. Finally, a certain isomorphism between Ore extensions are considered.
 }

\\

\end{minipage}
\end{center}
\normalsize

\section{Introduction}
\def\theequation{\thesection.\arabic{equation}}
\setcounter{equation}{0}

 A version of non-commutative polynomial ring, introduced by Ore in \cite{OO33}, has become one of the most basic and useful constructions in ring theory.  Such a polynomial ring is called Ore extension  now.
 From the perspective of quantum groups and Hopf algebras, the Ore extension is very important for constructing examples of Hopf algebras that are neither commutative nor cocommutative.
 These extensions are also called skew polynomial rings.

 In recent years, the Ore extension has been widely applied to other branches of Hopf algebra.
 Many new examples (usually finite dimensional) with special properties are constructed through Ore extension, such as pointed Hopf algebras, co-Frobenius Hopf algebras, and quasitriangular Hopf algebras.

 In \cite{PA03}, Panov, Aleksandr N. introduced the Hopf-Ore extensions, and gave the necessary and sufficient conditions for the Ore extension of a Hopf algebra to be a Hopf algebra.
 Li Chao and Li Jinqi \cite{LL20} introduced the concept of Ambikew Hopf $\pi$-coalgebra, which can be obtained from Hopf $\pi$-coalgebra through twice Ore extension.
 Jiao Zhengming \cite{JM14} further generalized the Hopf-Ore extensions theory to Hopf coquasigroups.
 After that, Wang Dingguo and Lu Daowei \cite{WL14} extended the Ore extension of Hopf algebras to the Hopf group coalgebras.
 Therefore, there is a natural question:  Does the Ore extension still hold for group-cograded Hopf coquasigroups? This is the motivation of our paper.

 For this question, we give a positive answer in this paper. The first matter we have to resolve is how to define the extension of group-cograded Hopf coquasigroups.
 This paper is organized as follows.

 In Section 2, we recall some concepts which will be used in the following section, such as Hopf coquasigroups, Ore extensions and the Turaev category.

 In Section 3, we provide the definition of group-cograded Hopf coquasigroups, and give a categorical interpretation and some interesting examples.

 In Section 4, we introduce the concept of the Ore extension for group-cograded Hopf coquasigroups, and give an equivalent condition that characterizes them as still being group-cograded Hopf coquasigroups. Finally, a isomorphism theorem of group-cograded Hopf coquasigroups is presented.

\section{Preliminaries}

 In this section, we review some basic definitions that need to be used in the following, such as Hopf coquasigroups, Ore extension of algebras, Ore extension of Hopf coquasigroups and Turaev category. Throughout this article, all spaces we considered are over a fixed field $k$.

\subsection{Hopf coquasigroups and Ore extension}

 First, let us recall the definition of Ore extension of an algebra from \cite{CK00}. Let $A$ be an algebra, and $\tau$ be an algebra endomorphism of $A$. A linear endomorphism $\delta$ of $A$ is called an $\tau$-derivation of $A$ if for all $a, b \in A$,
 \begin{eqnarray*}
 \delta(ab) = \delta(a)b+\tau(a)\delta(b).
 \end{eqnarray*}
 This condition implies $\delta (1) = 0$.

 The Ore extension $R = A[y; \tau, \delta]$ of $A$ is an algebra generated by the algebra $A$ and the variable $y$ with the relation
 \begin{eqnarray}
 ya=\tau(a)y + \delta(a), \label{2.3}
 \end{eqnarray}
 for all $a\in A$.
 \\

 Recall from \cite{KM10} that a Hopf coquasigroup $H$ is a unital associative algebra equipped with counital $\epsilon: H \rightarrow k$ and algebra homomorphisms $\Delta: H\rightarrow H\otimes H$ and linear map $S: H\rightarrow H$ such that
 \begin{eqnarray}
  &&(m\otimes id)(S\otimes id\otimes id)(id\otimes \Delta)\Delta
  = 1\otimes id
  = (m\otimes id)(id\otimes S\otimes id)(id\otimes \Delta)\Delta, \label{2.1}\\
  &&(id\otimes m)(id\otimes id\otimes S)(\Delta \otimes id)\Delta
  = id\otimes 1
  = (id\otimes m)(id\otimes S\otimes id)(\Delta \otimes id)\Delta. \label{2.2}
  \end{eqnarray}
 In this paper we use Sweedler notation, the conditions \eqref{2.1} and \eqref{2.2} come out as
 \begin{eqnarray*}
  && S(h_{(1)})h_{(2)(1)}\otimes h_{(2)(2)} = 1\otimes h
  = h_{(1)}S(h_{(2)(1)})\otimes h_{(2)(2)},\\
  && h_{(1)(1)}\otimes S(h_{(1)(2)})h_{(2)} = h\otimes 1
  = h_{(1)(1)}\otimes h_{(1)(2)}S(h_{(2)}),
 \end{eqnarray*}
 for all $h\in H$.
\\

 Let $H$ be a Hopf coquasigroup, $R = H[y; \tau, \delta]$ is called the Hopf coquasigroup-Ore extension in \cite{JM14} if $R$ is a Hopf coquasigroup with sub-Hopf coquasigroup $H$ and there exist $r_{1}, r_{2} \in H$ such that $\Delta(y) = y\otimes r_{1} + r_{2}\otimes y$.

\subsection{Turaev category}

 Turaev category as an special symmetric monoidal category is introduced by Caenepeel in \cite{CD06}. Let $K$ be a commutative ring.
 A Turaev $K$-module is a couple $\underline{M}=(X, M)$, where $X$ is a set,
 and $M=(M_{x})_{x\in X}$ is a family of $K$-modules indexed by $X$.
 A morphisms between two $T$-modules $(X, M)$ and $(Y, N)$ is a couple $\underline{\varphi}=(f,\varphi)$, where $f: Y\rightarrow X$ is a function,
 and $\varphi=(\varphi_{y}: M_{f_{(y)}}\rightarrow N_{y})_{y\in Y}$ is a family of linear maps indexed by $Y$.
 The composition of $\underline{\varphi}: \underline{M}\rightarrow \underline{N}$ and $\underline{\varphi}: \underline{N}\rightarrow \underline{P}=(Z,P)$ is defined as follows:
\begin{eqnarray*}
&&\underline{\psi}\circ \underline{\varphi}=(f\circ g, (\psi_{z}\circ \varphi_{g(z)})_{z\in Z}).
\end{eqnarray*}
The category of Turaev $K$-modules is called Turaev category and denoted by $\mathscr{T}_{K}$. Specially when $K$ is a field, then $\mathscr{T}_{K} = \mathscr{T}_{k}$.

\section{Group-cograded Hopf coquasigroups}

 Let $G$ be a group with the unit $1$. First of all, we introduce group-cograded Hopf coquasigroups.

 \textbf{Definition \thesection.1}
 $H = (\bigoplus_{p\in G}H_{p}, m, \mu, \Delta, \epsilon, S)$ is called a group-cograded Hopf coquasigroup over $k$, if the following conditions hold:
\begin{itemize}
  \item[(1)] each $H_{p}$ is an unital associative $k$-algebra with multiplication $m_{p}$ and unit $\mu_{p}$.

      $H_{p}H_{q}=0$ whenever $p,q\in G$ and $p\neq q$, and $\mu_{p}(1_{k})=1_{p}$;
  \item[(2)] comultiplication $\Delta$ is a family of homomorphisms $\{\Delta_{p, q}: H_{pq}\rightarrow H_{p}\otimes H_{q}\}_{p,q\in G}$,
      and counit $\epsilon: H_{1}\rightarrow k$ is an algeba homomorphism in the sense that for $p\in G,$
      \begin{eqnarray*}
      (id_{H_{p}}\otimes \epsilon)\Delta_{p,1} = (\epsilon\otimes id_{H_{p}})\Delta_{1,p} = id_{H_{p}},\\
      \epsilon(1_{1}) = 1_{k};
    \end{eqnarray*}
  \item[(3)] antipode $S$ is an algebra anti-homomorphism with $S = \{S_{p}: H_{p}\rightarrow H_{p^{-1}}\}_{p\in G}$, and for any $p, q\in G$,
    \begin{eqnarray}
      &&(m_{q}\otimes id_{H_{p}})(S_{q^{-1}}\otimes id_{H_{q}}\otimes id_{H_{p}})(id_{H_{q^{-1}}}\otimes \Delta_{q,p})\Delta_{q^{-1},qp}
      = \mu_{q}\otimes id_{H_{p}} \nonumber
      \\
      &=& (m_{q}\otimes id_{H_{p}})(id_{H_{q}}\otimes S_{q^{-1}}\otimes id_{H_{p}}) (id_{H_{q}}\otimes \Delta_{q^{-1},p})\Delta_{q,q^{-1}p}, \label{3.4}
      \\
      &&(id_{H_{p}} \otimes m_{q})(id_{H_{p}} \otimes id_{H_{q}} \otimes S_{q^{-1}})(\Delta_{p,q} \otimes id_{H_{q^{-1}}})\Delta_{pq,q^{-1}}
      = id_{H_{p}} \otimes \mu_{q}
      \nonumber
      \\
      &=& (id_{H_{p}} \otimes m_{q})(id_{H_{p}} \otimes S_{q^{-1}} \otimes id_{H_{q}})(\Delta_{p,q^{-1}} \otimes id_{H_{q}})\Delta_{pq^{-1},q}. \label{3.5}
    \end{eqnarray}
\end{itemize}

\textbf{Remark}
 (1) In the following, we use the Sweelder notation for the comultiplication:
 for any $p, q\in G$ and $h_{pq}\in H_{pq}$,
  \begin{eqnarray*}
    \Delta_{p,q}(h_{pq}) = h_{(1,p)}\otimes h_{(2,q)}.
  \end{eqnarray*}
 Then the conditions $\eqref{3.4}, \eqref{3.5}$ come out as
  \begin{eqnarray*}
    &&S_{q^{-1}}(h_{(1,q^{-1})})h_{(21,q)}\otimes h_{(22,p)} = 1_{q}\otimes h_{p}
    = h_{(1,q)}S_{q^{-1}}(h_{(21,q^{-1})})\otimes h_{(22,p)},
    \\
    &&h_{(11,p)}\otimes h_{(12,q)}S_{q^{-1}}(h_{(2,q^{-1})}) = h_{p}\otimes 1_{q}
    = h_{(11,p)}\otimes S_{q^{-1}}(h_{(12,q^{-1})})h_{(2,q)}.
  \end{eqnarray*}

  (2) If the comultiplication $\Delta$ of group-cograded Hopf coquasigroup $H$ is coassociative,
  then $H$ is actually a Hopf group-coalgebra introduced in \cite{V02}.

  (3) If $H=\bigoplus_{p\in G}H_{p}$ is a group-cograded Hopf coquasigroup with each component $H_{p}$ is finite dimensional, then $H^*=\bigoplus_{p\in G}H^*_{p}$ is a group-graded Hopf quasigroup introduced in \cite{SW18} with $Q$ be a group.
 \\

 As shown in \cite{CD06,LY22} Hopf group coalgbras(resp. group-graded Hopf quasigroups) are Hopf algebras(resp. Hopf quasigroups) in Turaev category $\mathcal T_k$.
 We give a categorical interpretation for group-cograded Hopf coquasigroups as follows.

\textbf{Proposition \thesection.2}
If $H=\bigoplus_{p\in G}H_{p}$ be a group-cograded Hopf coquasigroup, then $(G,H)$ is a Hopf coquasigroup in Turaev category $\mathcal T_k$.

\emph{Proof} As $H$ is a group-cograded Hopf coquasigroup and $G$ is a group, then we can give $\underline{H}=(G, H)$ a counital coalgebra structure $(\underline{H},\underline{\Delta},\underline{\epsilon})$ by
  \begin{eqnarray*}
    \begin{aligned}
      \underline{H}&\stackrel{\underline{\epsilon}}{\longrightarrow}\underline{k}
      \\
      G&\stackrel{i}{\longleftarrow} (*)
      \\
      H_{1}=H_{i(e)}&\stackrel{\epsilon}{\longrightarrow}k,
    \end{aligned}
    \qquad \rm{and} \qquad
    \begin{aligned}
      \underline{H}&\stackrel{\underline{\Delta}}{\longrightarrow}\underline{H}\otimes \underline{H}
      \\
      G&\stackrel{\eta}{\longleftarrow}G\times G
      \\
      H_{pq}=H_{\eta(p,q)}&\stackrel{\Delta_{p,q}}{\longrightarrow} H_{p}\otimes H_{q}.
    \end{aligned}
  \end{eqnarray*}

  We have
  \begin{eqnarray*}
    \begin{aligned}
      \underline{H}&\stackrel{\underline{\Delta}}{\makebox[1cm]{\rightarrowfill}}&\underline{H}\otimes \underline{H}&\stackrel{\underline{\epsilon}\otimes \underline{id}}{\makebox[1cm]{\rightarrowfill}}&\underline{H}
      \\
      G&\stackrel{\eta}{\makebox[1cm]{\leftarrowfill}}&G\times G&\stackrel{(i, G)}{\makebox[1cm]{\leftarrowfill}}&G
      \\
      H_{p}&\stackrel{\Delta_{e,p}}{\makebox[1cm]{\rightarrowfill}}&H_{1}\otimes H_{p}&\stackrel{\epsilon\otimes id}{\makebox[1cm]{\rightarrowfill}}&H_{p},
    \end{aligned}
    \qquad \rm{and} \qquad
    \begin{aligned}
      \underline{H}&\stackrel{\underline{\Delta}}{\makebox[1cm]{\rightarrowfill}}&\underline{H}\otimes \underline{H}&\stackrel{\underline{id}\otimes\underline{\epsilon}} {\makebox[1cm]{\rightarrowfill}}&\underline{H}
      \\
      G&\stackrel{\eta}{\makebox[1cm]{\leftarrowfill}}&G\times G &\stackrel{(G,i)}{\makebox[1cm]{\leftarrowfill}}&G
      \\
      H_{p}&\stackrel{\Delta_{p,e}}{\makebox[1cm]{\rightarrowfill}}&H_{p}\otimes H_{1}&\stackrel{id\otimes\epsilon}{\makebox[1cm]{\rightarrowfill}}&H_{p}.
    \end{aligned}
  \end{eqnarray*}

  Also, we can give $(G, H)$ a coalgebra structure $(\underline{H},\underline{\Delta},\underline{\epsilon})$ by
  \begin{eqnarray*}
    \begin{aligned}
      \underline{k}&\stackrel{\underline{\mu}}{\longrightarrow}\underline{H}
      \\
      (*)&\stackrel{e}{\longleftarrow} G
      \\
      k&\stackrel{\mu_{p}} {\longrightarrow}H_{p},
    \end{aligned}
    \qquad \rm{and} \qquad
    \begin{aligned}
      \underline{H}\otimes \underline{H}&\stackrel{\underline{m}}{\longrightarrow}\underline{H}
      \\
      G\times G&\stackrel{\delta}{\longleftarrow} G
      \\
      H_{p}\otimes H_{p}&\stackrel{\mu_{p}} {\longrightarrow}H_{p},
    \end{aligned}
  \end{eqnarray*}
  such that $(\underline{\Delta},\underline{\epsilon})$ are algebra maps.

  Let $s: G\rightarrow G, s(g)=g^{-1}$,  then we can consider a map $\underline{S}=(s,S)$ in Turaev category as the antipode of $\underline{H}$, where $S$ is the antipode of group-cograded Hopf coquasigroup $H$.
  The next we will check that $\underline{S}$ satisfy the condition (\ref{2.1}) and (\ref{2.2}) of Hopf quasigroup.
  \begin{eqnarray*}
    \begin{aligned}
      \underline{H}
      &\stackrel{\underline{\Delta}}{\makebox[1cm]{\rightarrowfill}}&
      \underline{H}\otimes \underline{H}
      &\stackrel {\underline{id}\otimes \underline{\Delta}}{\makebox[1cm]{\rightarrowfill}}& \underline{H}\otimes \underline{H}\otimes \underline{H}
      &\stackrel {\underline{S}\otimes \underline{id}\otimes \underline{id}}{\makebox[1cm]{\rightarrowfill}}&
      \underline{H}\otimes \underline{H}\otimes \underline{H}
      &\stackrel{\underline{m}\otimes \underline{id}}{\makebox[1cm]{\rightarrowfill}}& \underline{H}\otimes \underline{H}
      \\
      G
      &\stackrel{\eta}{\makebox[1cm]{\leftarrowfill}}&G\times G
      &\stackrel{(G,\eta)}{\makebox[1cm]{\leftarrowfill}}&G\times G\times G
      &\stackrel{(s,G,G)}{\makebox[1cm]{\leftarrowfill}}&G\times G\times G &\stackrel{(\delta,G)}{\makebox[1cm]{\leftarrowfill}}&G\times G
    \end{aligned}
  \end{eqnarray*}
  \begin{eqnarray*}
    \begin{aligned}
      H_{p}
      &\stackrel{\Delta_{q,q^{-1}p}}{\makebox[1cm]{\rightarrowfill}}&
      H_{q}\otimes H_{q^{-1}p}
      &\stackrel{id \otimes \Delta_{q^{-1},p}}{\makebox[1cm]{\rightarrowfill}}&
      H_{q}\otimes H_{q^{-1}}\otimes H_{p}
      &\stackrel{S_{q}\otimes id \otimes id}{\makebox[1cm]{\rightarrowfill}}&
      H_{q^{-1}}\otimes H_{q^{-1}}\otimes H_{p}  \\
      &\stackrel{m_{q^{-1}} \otimes id}{\makebox[1cm]{\rightarrowfill}}&
      H_{q^{-1}}\otimes H_{p}
    \end{aligned}
  \end{eqnarray*}
  and
  \begin{eqnarray*}
    \underline{k}\otimes \underline{H}&\stackrel{\underline{\mu}\otimes \underline{id}}{\makebox[1.5cm]{\rightarrowfill}}&\underline{H}\otimes \underline{H}
    \\
    (*)\times G&\stackrel{i\otimes G}{\makebox[1.5cm]{\leftarrowfill}}& G\times G
    \\
    k\otimes H_{p} &\stackrel{\mu_{p}\otimes id}{\makebox[1.5cm]{\rightarrowfill}}& H_{p}\otimes H_{p}.
  \end{eqnarray*}

  Since $H$ is a group-cograded Hopf coquasigroup, we have $(m\otimes id)(S\otimes id\otimes id)(id\otimes \Delta)\Delta = 1\otimes id$. Thus the left hand of equation (\ref{2.1}) holds, and the right hand is similar.

  The proof of equation (\ref{2.2}) for $\underline{H}$ is similar to the first one.
$\hfill \square$
\\

In the end of this section, we give some examples of group-cograded Hopf coquasigroups.

 \textbf{Example \thesection.3}
 Let $H$ be Hopf coquasigroup. Set $H^{G} = \bigoplus_{p\in G}H_{p}$ and $G$ is the homomorphism group of $H$, where for each $p \in G$, the algebra $H_{p}$ is a copy of $H$. Fix an identification isomorphism of algebras $i_{p}: H\rightarrow H_{p}$. For $p, q \in G$, we define a comultiplication
 $\Delta_{p, q}: H_{pq}\rightarrow H_{p}\otimes H_{q}$ by
 \begin{eqnarray*}
 &&\Delta_{p, q}(i_{pq}(h)) = i_{p}(h_{(1)})\otimes i_{q}(h_{(2)}),
 \end{eqnarray*}
 where $h\in H$. The counit $\epsilon: H_{e}\rightarrow k$ is defined by $\epsilon_{e}(i_{e}(h))=\epsilon(h)\in k$ for $h\in H$. For $p \in G$, the antipode $S_{p}: H_{p}\rightarrow H_{p^{-1}}$ is given by
 \begin{eqnarray*}
 && S_{p}(i_{p}(h))=i_{p^{-1}}(S(h)),
  \end{eqnarray*}
 where $h\in H$. It is easy to check that $H^{G}$ is a group-cograded Hopf coquasigroup.
\\

 \textbf{Example \thesection.4}
  Let $H=\bigoplus_{p\in G}H_{p}$ be a $Q$-graded Hopf quasigroup introduced in \cite{SW18}, where each component $H_{p}$ is finite dimensional and $Q$ is a group.
  Then $H^*=\bigoplus_{p\in G}H^*_{p}$ is a group-cograded Hopf coquasigroup.

\section{Ore extension}

The aim of this section is to prove a criterion for an Ore extension of a group-cograded Hopf coquasigroup to be a group-cograded Hopf coquasigroup with coproduct satisfying the relation $\eqref{3.6}$ bellow.
First we introduce the definition of Ore extension of a group-cograded Hopf coquasigroup.
\\

 \textbf{Definition \thesection.1}
 Let $H=\bigoplus_{p\in G}H_{p}$ be a group-cograded Hopf coquasigroup,
 the family $R = \bigoplus_{p\in G}R_{p} = \bigoplus_{p\in G}H_{p}[y_{p}; \tau_{p}, \delta_{p}]$ of $k$-spaces is called the group-cograded Hopf coquasigroup-Ore extension
 if $R = \bigoplus_{p\in G}R_{p}$ is also a group-cograded Hopf coquasigroup,
 where for any $p\in G$, $H_{p}[y_{p}; \tau_{p}, \delta_{p}]$ is the Ore extension of $H_{p}$,
 and there exist $r_{p}^{1}, r_{p}^{2}\in H_{p}$ such that
 \begin{eqnarray*}
  \Delta_{p,q}(y_{pq})=y_{p}\otimes r_{q}^{2} + r_{p}^{1}\otimes y_{q}.
 \end{eqnarray*}

 Note that $R_{1} = H_{1}[y_{1}; \tau_{1}, \delta_{1}]\}$ is the Hopf coquasigroup-Ore extension in the sense of \cite{JM14}.
 By definition 3.1 we have for $y_{p}\in H_{p}$,
 \begin{eqnarray*}
  && (m_{q}\otimes id_{H_{p}})(S_{q^{-1}}\otimes id_{H_{q}}\otimes id_{H_{p}})(id_{H_{q^{-1}}}\otimes \Delta_{q,p})\Delta_{q^{-1},qp}(y_{p})
  = 1_{q}\otimes y_{p} \\
  &=& (m_{q}\otimes id_{H_{p}})(id_{H_{q}}\otimes S_{q^{-1}}\otimes id_{H_{p}}) (id_{H_{q}}\otimes \Delta_{q^{-1},p})\Delta_{q,q^{-1}p}(y_{p}),  \\
  && (id_{H_{p}} \otimes m_{q})(id_{H_{p}} \otimes id_{H_{q}} \otimes S_{q^{-1}})(\Delta_{p,q} \otimes id_{H_{q^{-1}}})\Delta_{pq,q^{-1}}(y_{p})
  = y_{p} \otimes 1_{q} \\
  &=& (id_{H_{p}} \otimes m_{q})(id_{H_{p}} \otimes S_{q^{-1}} \otimes id_{H_{q}})(\Delta_{p,q^{-1}} \otimes id_{H_{q}})\Delta_{pq^{-1},q}(y_{p}).
\end{eqnarray*}
That is,
\begin{eqnarray*}
  &&S_{q^{-1}}(y_{q^{-1}})r_{(1,q)}^{2}\otimes r_{(2,p)}^{2} + S_{q^{-1}}(r_{q^{-1}}^{1})y_{q}\otimes r_{p}^{2} + S_{q^{-1}}(r_{q^{-1}}^{1})r_{q}^{1}\otimes y_{p} = 1_{q}\otimes y_{p}  \\
  &=& y_{q}S_{q^{-1}}(r_{(1,q^{-1})}^{2})\otimes r_{(2,p)}^{2} + r_{q}^{1}S_{q^{-1}}(y_{q^{-1}})\otimes r_{p}^{2} + r_{q}^{1}S_{q^{-1}}(r_{q^{-1}}^{1})\otimes y_{p},  \\
  &&y_{p}\otimes r_{q}^{2}S_{q^{-1}}(r_{q^{-1}}^{2}) + r_{p}^{1}\otimes y_{q}S_{q^{-1}}(r_{q^{-1}}^{2}) + r_{(1,p)}^{1}\otimes r_{(2,q)}^{1}S_{q^{-1}}(y_{q^{-1}}) = y_{p}\otimes 1_{q} \\
  &=& y_{p}\otimes S_{q^{-1}}(r_{q^{-1}}^{2})r_{q}^{2} + r_{p}^{1}\otimes S_{q^{-1}}(y_{q^{-1}})r_{q}^{2} + r_{(1,p)}^{1}\otimes S_{q^{-1}}(r_{(2,q^{-1})}^{1})y_{q}.
\end{eqnarray*}
So we have
\begin{eqnarray*}
  && S_{q^{-1}}(r_{q^{-1}}^{1})r_{q}^{1} = r_{q}^{1}S_{q^{-1}}(r_{q^{-1}}^{1}) = 1_{q}, \\
  && r_{q}^{2}S_{q^{-1}}(r_{q^{-1}}^{2}) = S_{q^{-1}}(r_{q^{-1}}^{2})r_{q}^{2} = 1_{q},
\end{eqnarray*}
and
\begin{eqnarray*}
&& S_{q^{-1}}(y_{q^{-1}})r_{(1,q)}^{2}\otimes r_{(2,p)}^{2}
     + S_{q^{-1}}(r_{q^{-1}}^{1})y_{q}\otimes r_{p}^{2} = 0\\
  &=& y_{q}S_{q^{-1}}(r_{(1,q^{-1})}^{2}) \otimes r_{(2,p)}^{2}
     + r_{q}^{1}S_{q^{-1}}(y_{q^{-1}})\otimes r_{p}^{2}, \\
  && r_{p}^{1}\otimes y_{q}S_{q^{-1}}(r_{q^{-1}}^{2})
     + r_{(1,p)}^{1}\otimes r_{(2,q)}^{1}S_{q^{-1}}(y_{q^{-1}}) = 0  \\
  &=& r_{p}^{1}\otimes S_{q^{-1}}(y_{q^{-1}})r_{q}^{2}
     + r_{(1,p)}^{1}\otimes S_{q^{-1}}(r_{(2,q^{-1})}^{1})y_{q}.
\end{eqnarray*}

As in Section 3 of \cite{JM14} and Lemma 3.3 of \cite{WL14}, if $r^{1}, r^{2}$ satisfies the following conditions
\begin{eqnarray*}
  \Delta_{p,q}(r_{pq}^{1}) = r_{p}^{1}\otimes r_{q}^{1}, \quad
  \Delta_{p,q}(r_{pq}^{2}) = r_{p}^{2}\otimes r_{q}^{2},
\end{eqnarray*}
 then we have
\begin{eqnarray*}
  (r_{q}^{i})^{-1}
  &=& S_{q^{-1}}(r_{q^{-1}}^{i}), \quad i=1,2,
  \\
  (S_{q^{-1}}(y_{q^{-1}})r_{q}^{2} + S_{q^{-1}}(r_{q^{-1}}^{1})y_{q})\otimes r_{p}^{2}
  &=& (y_{q}S_{q^{-1}}(r_{q^{-1}}^{2}) + r_{q}^{1}S_{q^{-1}}(y_{q^{-1}}))\otimes r_{p}^{2} = 0,
  \\
  r_{p}^{1}\otimes (y_{q}S_{q^{-1}}(r_{q^{-1}}^{2}) + r_{q}^{1}S_{q^{-1}}(y_{q^{-1}}))
  &=& r_{p}^{1}\otimes (S_{q^{-1}}(y_{q^{-1}})r_{q}^{2} + S_{q^{-1}}(r_{q^{-1}}^{1})y_{q}) = 0.
\end{eqnarray*}
Replacing the generating elements $y_{p}$ by $y_{p}'=y_{p}(r_{p}^{2})^{-1}$ and $r_{p}^{1}(r_{p}^{2})^{-1}$ by $r_{p}$, we see that
\begin{eqnarray*}
  \Delta_{p,q}(y_{pq}') &=& \Delta_{p,q}(y_{pq}(r_{pq}^{2})^{-1})
  \\
  &=& \Delta_{p,q}(y_{pq})\Delta_{p,q}(r_{pq}^{2})^{-1})
  \\
  &=& (y_{p}\otimes r_{q}^{2} + r_{p}^{1}\otimes y_{q})((r_{p}^{2})^{-1}\otimes (r_{q}^{2})^{-1})
  \\
  &=& y_{p}(r_{p}^{2})^{-1} \otimes r_{q}^{2}(r_{q}^{2})^{-1} + r_{p}^{1}(r_{p}^{2})^{-1}\otimes y_{q}(r_{q}^{2})^{-1}
  \\
  &=& y_{p}' \otimes 1_{q} + r_{p} \otimes y_{q}'.
\end{eqnarray*}
So we always assume in what follows that the elements $y_{p}$ in the group-cograded Hopf coquasigroup-Ore extension satisfying the relations
\begin{eqnarray}
  \Delta_{p,q}(y_{pq}) = y_{p} \otimes 1_{q} + r_{p} \otimes y_{q},  \label{3.6}
\end{eqnarray}
for some elements $r_{p} \in H_{p}$. As usual, $Ad_{r_{p}}(h) = r_{p}hS_{p^{-1}}(r_{p^{-1}}) = r_{p}h(r_{p})^{-1}$.
\\

\textbf{Lemma \thesection.2}
Let $H=\bigoplus_{p\in G}H_{p}$ be a group-cograded Hopf coquasigroup. If $R = \bigoplus_{p\in G}R_{p} = \bigoplus_{p\in G}H_{p}[y_{p}; \tau_{p}, \delta_{p}]$ is the group-cograded Hopf coquasigroup-Ore extension of $H$, then
\begin{eqnarray}
  S_{p^{-1}}(y_{p^{-1}}) = -(r_{p})^{-1}y_{p},  \label{3.7}
\end{eqnarray}
where $(r_{p})^{-1} = S_{p^{-1}}(r_{p^{-1}})$.
\\

\emph{Proof}
By \eqref{3.6}, we have
\begin{eqnarray*}
m_{p}(S_{p^{-1}} \otimes id_{H_{p}}) \Delta_{p^{-1},p}(y_{1}) = \epsilon(y_{1})1_{p} = 0,
\end{eqnarray*}
then
\begin{eqnarray*}
S_{p^{-1}}(y_{p^{-1}}) + S_{p^{-1}}(r_{p^{-1}})y_{p} = 0.
\end{eqnarray*}
So
$S_{p^{-1}}(y_{p^{-1}}) = - S_{p^{-1}}(r_{p^{-1}})y_{p} = -(r_{p})^{-1}y_{p}$.
$\hfill \square$
\\

 Following the above results, we obtain the main theorem of this paper.

\textbf{Theorem \thesection.3}
Let $H=\bigoplus_{p\in G}H_{p}$ be a group-cograded Hopf coquasigroup. Then the group-cograded Hopf coquasigroup $R = \bigoplus_{p\in G}R_{p} = \bigoplus_{p\in G}H_{p}[y_{p}; \tau_{p}, \delta_{p}]$ is the group-cograded Hopf coquasigroup-Ore extension if and only if

\begin{itemize}
  \item[(1)] there is a character $\chi : H_{1}\rightarrow k$ such that for any $p\in G$
    \begin{eqnarray}
      \tau_{p}(h_{p}) = \chi(h_{(1,1)})h_{(2,p)}, \label{D1}
    \end{eqnarray}
    where $h_{p}\in H_{p}$;
  \item[(2)] the following relations hold:
    \begin{eqnarray}
      \chi(h_{(1,1)})h_{(21,p)}\otimes h_{(22.q)} &=& Ad_{r_{p}}(h_{(1,p)})\chi(h_{(21,1)}) \otimes h_{(22,q)}  \nonumber \\
      &=& \chi(h_{(11,1)})h_{(12,p)}\otimes h_{(2,q)}, \label{D2}
    \end{eqnarray}
  \item[(3)] the $\tau_{p}$-derination $\delta_{p}$ satisfies the relation
    \begin{eqnarray}
      \Delta_{p,q}(\delta_{pq}(h_{p,q}))=\delta_{p}(h_{(1,p)})\otimes h_{(2,q)} + r_{p}h_{(1,p)}\otimes \delta_{q}(h_{(2,q)}). \label{D3}
    \end{eqnarray}
\end{itemize}

\emph{Proof}
 The proof is divided into three parts.
 At step 1 we show that the comultiplication $\Delta = \bigoplus_{p\in G}\Delta_{p}$ can be extended to $R = \bigoplus_{p\in G}R_{p} = \bigoplus_{p\in G}H_{p}[y_{p}; \tau_{p}, \delta_{p}]$ by \eqref{3.6} if and only if relations $\eqref{D1}-\eqref{D3}$ hold.
 At step 2 we prove that $R_{1}$ admits an extension of the counit from $H_{1}$(in fact this has been proved in [6]).
 At step 3 we show that $R$ has antipode $S$ extending the antipode $S|_{H}$ by \eqref{3.7}.

  \begin{itemize}
    \item[Step 1.]  $Comultiplication.$

    Assume that the comultiplication $\Delta |_{H}$ can be extended to $R = \bigoplus_{p\in G}R_{p} = \bigoplus_{p\in G}H_{p}[y_{p}; \tau_{p}, \delta_{p}]$ \eqref{3.6}. Then the homomorphism $\Delta$ preserve the relation
    \begin{eqnarray}
      y_{p}h_{p} = \tau_{p}(h_p)y_{p} + \delta_{p}(h_p), \label{E1}
    \end{eqnarray}
    for any $p\in G$ and $h_p\in H_{p}$, i.e.,
    \begin{eqnarray}
      \Delta_{p,q}(y_{pq})\Delta_{p,q}(h_{pq}) = \Delta_{p,q}(\tau_{pq}(h_{pq}))\Delta_{p,q}(y_{pq}) + \Delta_{p,q}(\delta_{pq}(h_{pq})), \label{E2}
    \end{eqnarray}
    for any $h_{pq}\in H_{pq}$. We have
    \begin{eqnarray*}
      \Delta_{p,q}(y_{pq})\Delta_{p,q}(h_{pq})
      &=& (y_{p} \otimes 1_{q} + r_{p} \otimes y_{q})(h_{(1,p)}\otimes h_{(2,q)})
      \\
      &=& y_{p}h_{(1,p)}\otimes h_{(2,q)} + r_{p}h_{(1,p)}\otimes y_{q}h_{(2,q)}
      \\
      &=& \tau_{p}(h_{(1,p)})y_{p}\otimes h_{(2,q)} + \delta_{p}(h_{(1,p)})\otimes h_{(2,q)}
      \\
      && + r_{p}h_{(1,p)}\otimes \tau_{q}(h_{(2,q)})y_{q} + r_{p}h_{(1,p)}\otimes \delta_{q}(h_{(2,q)})
      \\
      &=& (\tau_{p}(h_{(1,p)})\otimes h_{(2,q)})(y_{p}\otimes 1_{q}) \\
      &&+ (r_{p}h_{(1,p)}r_{p}^{-1}\otimes \tau_{q}(h_{(2,q)}))(r_{p}\otimes y_{q})
      \\
      && + \delta_{p}(h_{(1,p)})\otimes h_{(2,q)} + r_{p}h_{(1,p)}\otimes \delta_{q}(h_{(2,q)}),
    \end{eqnarray*}
    and
    \begin{eqnarray*}
      &&\Delta_{p,q}(\tau_{pq}(h_{pq}))\Delta_{p,q}(y_{pq}) + \Delta_{p,q}(\delta_{pq}(h_{pq}))
      \\
      =&&\Delta_{p,q}(\tau_{pq}(h_{pq}))(y_{p} \otimes 1_{q} + r_{p} \otimes y_{q}) + \Delta_{p,q}(\delta_{pq}(h_{pq}))
      \\
      =&&\Delta_{p,q}(\tau_{pq}(h_{pq}))(y_{p} \otimes 1_{q} ) + \Delta_{p,q}(\tau_{pq}(h_{pq}))(r_{p} \otimes y_{q}) + \Delta_{p,q}(\delta_{pq}(h_{pq})).
    \end{eqnarray*}
    It is clear for $\Delta$ to preserves \eqref{E1} if and only if the following relations hold:
    \begin{eqnarray}
      \Delta_{p,q}(\tau_{pq}(h_{pq})) &=& \tau_{p}(h_{(1,p)})\otimes h_{(2,q)}, \label{F1}
      \\
      \Delta_{p,q}(\tau_{pq}(h_{pq})) &=& Ad_{r_{p}}(h_{(1,p)})\otimes \tau_{q}(h_{(2,q)}), \label{F2}
      \\
      \Delta_{p,q}(\delta_{pq}(h_{pq})) &=& \delta_{p}(h_{(1,p)})\otimes h_{(2,q)} + r_{p}h_{(1,p)}\otimes \delta_{q}(h_{(2,q)}), \nonumber
    \end{eqnarray}
    for any $h\in H_{pq}$.The last equation coincides with \eqref{D3}.

    Now we shou that \eqref{F1} and \eqref{F2} imply \eqref{D1} and \eqref{D2}. Put $\chi(h_{1}):=\epsilon(\tau_{1}(h_{1}))$. It is clear that $\chi(h_{1})\in k$. One can regard $\chi$ as a mapping $\chi : H_1\rightarrow k$. Since for any $p\in G$, $\tau$ is an endomorphism, it follows that
    \begin{eqnarray*}
      \chi (h_{1}g_{1}) &=& \epsilon (\tau_{1} (h_{1}g_{1})) = \epsilon (\tau_{1}(h_{1})) \epsilon (\tau_{1}(g_{1})) =  \chi (h_{1}) \chi (g_{1}),
      \\
      \chi (h_{1} + g_{1}) &=& \epsilon (\tau_{1} (h_{1} + g_{1})) = \epsilon (\tau_{1} (h_{1})) + \epsilon (\tau_{1} (g_{1})) = \chi (h_{1}) + \chi (g_{1}).
    \end{eqnarray*}
    For $\tau_{p}(h_{p})\in H_{p}$ by \eqref{F1}, we have
    \begin{eqnarray*}
      (id_{H_{1}}\otimes \Delta_{1,p})\Delta_{1,p}(\tau_{p}(h_{p})) = \tau_{1}(h_{(1,1)}) \otimes h_{(21,1)} \otimes h_{(22,p)}.
    \end{eqnarray*}
    Then by \eqref{3.4} we have
    \begin{eqnarray}
      S_{1}(\tau_{1}(h_{(1,1)}))(h_{(21,1)}) \otimes h_{(22,p)} = 1\otimes \tau_{p}(h_{p}) = \tau_{1}(h_{(1,1)}) S_{1}( (h_{(21,1)})) \otimes h_{(22,p)}. \label{3.15}
    \end{eqnarray}
    Using the formula above, one can recover $\tau$ from $\chi$,
    \begin{eqnarray*}
      \chi(h_{(1,1)})h_{(2,p)}
      &=& \epsilon(\tau_{1}(h_{(1,1)}))h_{(2,p)}
      \\
      &=& \epsilon(\tau_{1}(h_{(1,1)}))\epsilon(h_{(21,1)})h_{(22,p)}
      \\
      &=& \epsilon(\tau_{1}(h_{(1,1)}))\epsilon (S_{1}(h_{(21,1)}))h_{(22,p)}
      \\
      &=& \epsilon(\tau_{1}(h_{(1,1)})S_{1}(h_{(21,1)}))h_{(22,p)}
      \\
      &\stackrel{\eqref{3.15}}=& \epsilon(1_{1})\tau_{p}(h_{p})
      \\
      &=& \tau_{p}(h_{p}).
    \end{eqnarray*}
    This proves \eqref{D1}.
    Substituting $\tau_{p}(h_{p})$ into \eqref{F1}, we obtain
    \begin{eqnarray*}
      \Delta_{p,q}(\tau_{pq}(h_{pq}))
      &=& \tau_{p}(h_{(1,p)})\otimes h_{(2,q)}
      \\
      &\stackrel{\eqref{D1}}=& \chi(h_{11,1})h_{(12,p)}\otimes h_{(2,q)},
    \end{eqnarray*}
    and substituting $\tau_{p}(h_{p})$ into \eqref{F2}, we obtain
    \begin{eqnarray*}
      \Delta_{p,q}(\tau_{pq}(h_{pq}))
      &\stackrel{\eqref{D1}}=& \Delta_{p,q}(\chi(h_{(1,1)})h_{(2,pq)})
      \\
      &=& \chi(h_{(1,1)}h_{(21,p)} \otimes h_{(22,q)}
      \\
      &\stackrel{\eqref{F2}}=& Ad_{r_{p}}(h_{(1,p)}) \otimes \tau_{q}(h_{(2,q)})
      \\
      &\stackrel{\eqref{D1}}=& Ad_{r_{p}}(h_{(1,p)}) \otimes \chi(h_{(21,1)})h_{(22,q)}.
    \end{eqnarray*}
    Hence
    \begin{eqnarray*}
      \chi(h_{(1,1)})h_{(21,p)}\otimes h_{(22.q)} = Ad_{r_{p}}(h_{(1,p)})\chi(h_{(21,1)}) \otimes h_{(22,q)} = \chi(h_{(11,1)})h_{(12,p)}\otimes h_{(2,q)}.
    \end{eqnarray*}
    This proves \eqref{D2}.We have proved that conditions $\eqref{D1}-\eqref{D3}$ are necessary
    conditions of the comultiplication.

    On the other hand, if conditions $\eqref{D1}-\eqref{D3}$ hold, then
    \begin{eqnarray*}
      \Delta_{p,q}(\tau_{pq}(h_{pq}))
      &\stackrel{\eqref{D1}}=& \Delta_{p,q}(\chi(h_{(1,1)})h_{(2,pq)})
      \\
      &=& \chi(h_{(1,1)})h_{(21,p)} \otimes h_{(22,q)}
      \\
      &\stackrel{\eqref{D2}}=& \chi(h_{(11,1)})h_{(12,p)} \otimes h_{(2,q)}
      \\
      &\stackrel{\eqref{D1}}=& \tau_{p}(h_{(1,p)}) \otimes h_{(2,q)},
      \\
      \Delta_{p,q}(\tau_{pq}(h_{pq}))
      &\stackrel{\eqref{D1}}=& \Delta_{p,q}(\chi(h_{(1,1)})h_{(2,pq)})
      \\
      &=& \chi(h_{(1,1)})h_{(21,p)} \otimes h_{(22,q)}
      \\
      &\stackrel{\eqref{D2}}=& Ad_{r_{p}}(h_{(1,p)})\chi(h_{(21,1)}) \otimes h_{(22,q)}
      \\
      &\stackrel{\eqref{D1}}=& Ad_{r_{p}}(h_{(1,p)}) \otimes \tau_{q}(h_{(2,q)}).
    \end{eqnarray*}
    This proves the relations \eqref{F1} and \eqref{F2} hold and the comultiplication $\Delta|_{H}$ can be extended to a homomorphism $\Delta: R\rightarrow R\otimes R$, and $\Delta$ is not required to be coassociative.

    \item[Step 2.]  $Counit.$

     For this part, from \cite{JM14} we have known that, as $R_{1}$ admits a comultiplication, there exists a counit extending $\epsilon|_{H_{1}}$ and satisfying $\epsilon(y_{1}) = 0$. It follows that $\epsilon$ admits an extension to R if and only if
     \begin{eqnarray*}
       \epsilon(\delta_{1}(h_{1})) = 0,
     \end{eqnarray*}
     for any $h_{1}\in H_{1}$.

    \item[Step 3.]  $Antipode.$

    Let $R$ be as in Step 1. Recall that $S = \bigoplus_{p\in G}S_{p}: H_{p}\rightarrow H_{p^{-1}}$ with $S_{p}$ being an antiautomorphism. If $R$ admits an antipode $S$ which can be extended from $H$ to $R$ by means of \eqref{3.7}, then $S$ satisfies \eqref{3.4} and \eqref{3.5} and preserve \eqref{E1}. This means that for any $h\in H_{p}$
    \begin{eqnarray}
      && (m_{q}\otimes id_{H_{p}})(S_{q^{-1}}\otimes id_{H_{q}}\otimes id_{H_{p}})(id_{H_{q^{-1}}}\otimes \Delta_{q,p})\Delta_{q^{-1},qp}(y_{p}h)
      = 1_{q}\otimes y_{p}h
      \nonumber
      \\
      && = (m_{q}\otimes id_{H_{p}})(id_{H_{q}}\otimes S_{q^{-1}} \otimes id_{H_{p}}) (id_{H_{q}}\otimes \Delta_{q^{-1},p})\Delta_{q,q^{-1}p}(y_{p}h), \label{3.16}
      \\
      && (id_{H_{p}} \otimes m_{q})(id_{H_{p}} \otimes id_{H_{q}} \otimes S_{q^{-1}})(\Delta_{p,q} \otimes id_{H_{q^{-1}}})\Delta_{pq,q^{-1}}(y_{p}h)
       = y_{p}h \otimes 1_{q} \nonumber
       \\
      && = (id_{H_{p}} \otimes m_{q})(id_{H_{p}} \otimes S_{q^{-1}} \otimes id_{H_{q}})(\Delta_{p,q^{-1}} \otimes id_{H_{q}})\Delta_{pq^{-1},q}(y_{p}h), \label{3.17}
    \end{eqnarray}
    and
    \begin{eqnarray}
      S_{p}(h)S_{p}(y_{p}) = S_{p}(y_{p})S_{p}(\tau_{p}(h)) + S_{p}(\delta_{p}(h)).  \label{3.18}
    \end{eqnarray}

    Let us prove $\eqref{3.16}-\eqref{3.17}$ :
    \begin{eqnarray*}
      && (m_{q}\otimes id_{H_{p}})(S_{q^{-1}}\otimes id_{H_{q}}\otimes id_{H_{p}})(id_{H_{q^{-1}}}\otimes \Delta_{q,p})\Delta_{q^{-1},qp}(y_{p}h)
      \\
      &\stackrel{\eqref{E1}}=& (m_{q}\otimes id_{H_{p}})(S_{q^{-1}}\otimes id_{H_{q}}\otimes id_{H_{p}})(id_{H_{q^{-1}}}\otimes \Delta_{q,p}) \\
      &&\Delta_{q^{-1},qp}(\tau_{p}(h)y_{p} + \delta_{p}(h))
      \\
      &\stackrel{\eqref{3.6}\eqref{F1}}=&(m_{q}\otimes id_{H_{p}})(S_{q^{-1}}\otimes id_{H_{q}}\otimes id_{H_{p}})(id_{H_{q^{-1}}}\otimes \Delta_{q,p})\\
      &&(\tau_{q^{-1}}(h_{(1,q^{-1})})y_{q^{-1}}\otimes h_{(2,qp)} + \tau_{q^{-1}}(h_{(1,q^{-1})})r_{q^{-1}}\otimes h_{(2,qp)}y_{qp} + \Delta_{q^{-1},qp}(\delta_{p}(h)))
      \\
      &=& (m_{q}\otimes id_{H_{p}})(S_{q^{-1}}\otimes id_{H_{q}}\otimes id_{H_{p}})(id_{H_{q^{-1}}} \otimes \Delta_{q,p}) \\
      &&(\tau_{q^{-1}}(h_{(1,q^{-1})})y_{q^{-1}}\otimes h_{(2,qp)} + \tau_{q^{-1}}(h_{(1,q^{-1})})r_{q^{-1}}\otimes h_{(2,qp)}y_{qp})\\
      && + (m_{q}\otimes id_{H_{p}})(S_{q^{-1}}\otimes id_{H_{q}}\otimes id_{H_{p}})(id_{H_{q^{-1}}} \otimes \Delta_{q,p})\Delta_{q^{-1},qp}(\delta_{p}(h))
      \\
      &\stackrel{\eqref{3.4}}=& S_{q^{-1}}(y_{q^{-1}})S_{q^{-1}}(\tau_{q^{-1}}(h_{(1,q^{-1})}))h_{(21,q)} \otimes h_{(22,p)} \\
      &&+ S_{q^{-1}}(r_{q^{-1}})S_{q^{-1}}(\tau_{q^{-1}}(h_{(1,q^{-1})}))h_{(21,q)}y_{q} \otimes h_{(22,p)} \\
      &&+ S_{q^{-1}}(r_{q^{-1}})S_{q^{-1}}(\tau_{q^{-1}}(h_{(1,q^{-1})}))h_{(21,q)}r_{q} \otimes h_{(22,p)}y_{p} + 1_{q}\otimes \delta_{p}(h)
      \\
      &\stackrel{\eqref{3.15}}=& S_{q^{-1}}(y_{q^{-1}})\otimes \tau_{p}(h) + S_{q^{-1}}(r_{q^{-1}})y_{q}\otimes \tau_{p}(h) + S_{q^{-1}}(r_{q^{-1}})r_{q}\otimes \tau_{p}(h)y_{p} \\
      &&+ 1_{q}\otimes \delta_{p}(h)
      \\
      &=& S_{q^{-1}}(y_{q^{-1}})\otimes \tau_{p}(h) + S_{q^{-1}}(r_{q^{-1}})y_{q}\otimes \tau_{p}(h) + 1\otimes \tau_{p}(h)y_{p} + 1_{q}\otimes \delta_{p}(h)
      \\
      &=& (S_{q^{-1}}(y_{q^{-1}})+ S_{q^{-1}}(r_{q^{-1}})y_{q})\otimes \tau_{p}(h) + 1_{q}\otimes (\tau_{p}(h)y_{p} + \delta_{p}(h))
      \\
      &=& 1_{q}\otimes y_{p}h.
    \end{eqnarray*}
    This proves \eqref{3.16}. Relation \eqref{3.17} can be proved similarly.

    For \eqref{3.18}, it follows from \eqref{3.7} that for any $h\in H_{p}$
    \begin{eqnarray*}
      -S_{p}(h)(r_{p^{-1}})^{-1}y_{p^{-1}} &=& -(r_{p^{-1}})^{-1}y_{p^{-1}}S_{p}(\tau_{p}(h)) + S_{p}(\delta_{p}(h)),
      \\
      -S_{p}(h)(r_{p^{-1}})^{-1}y_{p^{-1}} &=& -(r_{p^{-1}})^{-1}\tau_{p^{-1}}(S_{p}(\tau_{p}(h)))y_{p^{-1}} \\
      &&- (r_{p^{-1}})^{-1}\delta_{p^{-1}}(S_{p}(\tau_{p}(h))) + S_{p}(\delta_{p}(h)).
    \end{eqnarray*}
    Condition \eqref{3.18} holds if and only if the following two conditions hold:
    \begin{eqnarray}
      S_{p}(h)(r_{p^{-1}})^{-1} &=& (r_{p^{-1}})^{-1}\tau_{p^{-1}}(S_{p}(\tau_{p}(h))), \label{3.19}  \\
      r_{p^{-1}}S_{p}(\delta_{p}(h)) &=& \delta_{p^{-1}}(S_{p}(\tau_{p}(h))). \label{3.20}
    \end{eqnarray}
    Let us prove \eqref{3.19}. We have
    \begin{eqnarray*}
      &&\tau_{p^{-1}}(S_{p}(\tau_{p}(h)))
      \\
      &\stackrel{\eqref{D1}}=& \tau_{p^{-1}}(S_{p}(\chi(h_{(1,1)})h_{(2,p)}))
      \\
      &=& \chi(h_{(1,1)})\tau_{p^{-1}}(S_{p}(h_{(2,p)}))
      \\
      &\stackrel{\eqref{D1}}=& \chi(h_{(1,1)})\chi(S_{p}(h_{(2,p)})_{(1,1)})S_{p}(h_{(2,p)})_{(2,p^{-1})}
      \\
      &=& \chi(h_{(1,1)})\chi(S_{1}(h_{(22,1)}))S_{p}(h_{(21,p)})
      \\
      &=& \chi(S_{1}(h_{(22,1)}))S_{p}(\chi(h_{(1,1)})h_{(21,p)})
      \\
      &\stackrel{\eqref{D2}}=& \chi(S_{1}(h_{(22,1)}))S_{p}(Ad_{r_{p}}(h_{(1,p)})\chi(h_{(21,1)}))
      \\
      &=& S_{p}(Ad_{r_{p}}(h_{(1,p)}))\chi(h_{(21,1)}S_{1}(h_{(22,1)}))
      \\
      &=& S_{p}(Ad_{r_{p}}(h))
      \\
      &=& r_{p^{-1}}S_{p}(h)(r_{p^{-1}})^{-1}
      \\
      &=& Ad_{r_{p^{-1}}}(S_{p}(h)).
    \end{eqnarray*}
    Our next objective is to prove \eqref{3.20}. It follows from \eqref{D1} that
    \begin{eqnarray*}
      S_{p}(\tau_p(h)) = S_{p}(\chi(h_{(1,1)})h_{(2,p)}).
    \end{eqnarray*}
    Therefore, \eqref{3.20} can be represented in equivalent form as
    \begin{eqnarray}
      r_{p^{-1}}S_{p}(\delta_{p}(h)) = \chi(h_{(1,1)})\delta_{p^{-1}}(S_{p}(h_{(2,p)})). \label{3.21}
    \end{eqnarray}
    We denote $L_{p} = r_{p^{-1}}S_{p}(\delta_{p}(h))$ and $M_{p} = \chi(h_{(1,1)})\delta_{p^{-1}}(S_{p}(h_{(2,p)}))$.

    From \eqref{D3} we have
    \begin{eqnarray*}
      \Delta_{p,p^{-1}}(\delta_{1}(h_{1}))=\delta_{p}(h_{(1,p)})\otimes h_{(2,p^{-1})} + r_{p}h_{(1,p)}\otimes \delta_{p^{-1}}(h_{(2,p^{-1})}),
    \end{eqnarray*}
    and we apply $m_{p}(id_{H_{p}}\otimes S_{p^{-1}})$ to the above equality, we obtain
    \begin{eqnarray*}
      m_{p}(id_{H_{p}}\otimes S_{p^{-1}})(\Delta_{p,p^{-1}}(\delta_{1}(h_{1})))
      &=& m_{p}(id_{H_{p}}\otimes S_{p^{-1}})(\delta_{p}(h_{(1,p)})\otimes h_{(2,p^{-1})} \\
      &&+ r_{p}h_{(1,p)}\otimes \delta_{p^{-1}}(h_{(2,p^{-1})})),
      \\
      0 = \epsilon(\delta_{1}(h_{1}))1_{p} &=& \delta_{p}(h_{(1,p)})S_{p^{-1}}(h_{(2,p^{-1})}) \\
      &&+ r_{p}h_{(1,p)}S_{p^{-1}}(\delta_{p^{-1}}(h_{(2,p^{-1})})).
    \end{eqnarray*}
    Thus
    \begin{eqnarray}
      -r_{p}^{-1}\delta_{p}(h_{(1,p)})S_{p^{-1}}(h_{(2,p^{-1})}) = h_{(1,p)}S_{p^{-1}}(\delta_{p^{-1}}(h_{(2,p^{-1})})). \label{3.22}
    \end{eqnarray}
    Then for any $h_{p^{-1}}\in H_{p^{-1}}$
    \begin{eqnarray*}
      L_{p^{-1}} &=& r_{p}S_{p^{-1}}(\delta_{p^{-1}}(h_{p^{-1}}))
      \\
      &
      =& r_{p}S_{p^{-1}}(h_{(1,p^{-1})})h_{(21,p)}S_{p^{-1}}(\delta_{p^{-1}}(h_{(22,p^{-1})}))
      \\
      &\stackrel{\eqref{3.22}}=& -r_{p}S_{p^{-1}}(h_{(1,p^{-1})})(r_{p})^{-1}\delta_{p}(h_{(21,p)})S_{p^{-1}}(h_{(22,p^{-1})})
      \\
      &=& -Ad_{r_{p}}(S_{p^{-1}}(h_{(1,p^{-1})}))\delta_{p}(h_{(21,p)})S_{p^{-1}}(h_{(22,p^{-1})}).
    \end{eqnarray*}
    On the other hand, for any $h_{1}\in H_{1}$, we have $\epsilon(h_{1})1_{p} = h_{(1,p)}S_{p^{-1}}(h_{(2,p^{-1})})$. The action by $\delta_{p}$ on both sides gives
    \begin{eqnarray}
      0 = \delta_{p}(h_{(1,p)})S_{p^{-1}}(h_{(2,p^{-1})}) + \tau_{p}(h_{(1,p)})\delta_{p}(S_{p^{-1}}(h_{(2,p^{-1})})).\label{3.23}
    \end{eqnarray}
    Then we have
    \begin{eqnarray*}
      M_{p^{-1}} &=& \chi(h_{(1,1)})\delta_{p}(S_{p^{-1}}(h_{(2,p^{-1})}))
      \\
      &=& \chi(h_{(1,1)})\tau_{p}(S_{p^{-1}}(h_{(21,{p^{-1}})})h_{(221,p)}) \delta_{p}(S_{p^{-1}}(h_{(222,p^{-1})}))
      \\
      &=&\chi(h_{(1,1)})\tau_{p}(S_{p^{-1}}(h_{(21,{p^{-1}})}))\tau_{p}(h_{(221,p)}) \delta_{p}(S_{p^{-1}}(h_{(222,p^{-1})}))
      \\
      &\stackrel{\eqref{3.23}}=& -\chi(h_{(1,1)})\tau_{p}(S_{p^{-1}}(h_{(21,{p^{-1}})})) \delta_{p}(h_{(221,p)})S_{p^{-1}}(h_{(222,{p^{-1}})})
      \\
      &=& -\tau_{p}(S_{p^{-1}}(\chi(h_{(1,1)})h_{(21,{p^{-1}})})) \delta_{p}(h_{(221,p)})S_{p^{-1}}(h_{(222,{p^{-1}})})
      \\
      &=& -\tau_{p}(S_{p^{-1}}(\chi(h_{(11,1)})h_{(12,{p^{-1}})})) \delta_{p}(h_{(21,p)})S_{p^{-1}}(h_{(22,{p^{-1}})})
      \\
      &=& -\chi(h_{(11,1)})\tau_{p}(S_{p^{-1}}(h_{(12,{p^{-1}})})) \delta_{p}(h_{(21,p)})S_{p^{-1}}(h_{(22,{p^{-1}})})
      \\
      &=& -\chi(h_{(11,1)})Ad_{r_{p}}(S_{p^{-1}}(h_{(122,{p^{-1}})}))\chi(S_{1}(h_{(121,1)})) \delta_{p}(h_{(21,p)}) \\
      &&S_{p^{-1}}(h_{(22,{p^{-1}})})
      \\
      &=& -\chi(h_{(11,1)}S_{1}(h_{(121,1)}))Ad_{r_{p}}(S_{p^{-1}}(h_{(122,{p^{-1}})})) \delta_{p}(h_{(21,p)})S_{p^{-1}}(h_{(22,{p^{-1}})})
      \\
      &=& -Ad_{r_{p}}(S_{p^{-1}}(h_{(1,p^{-1})}))\delta_{p}(h_{(21,p)})S_{p^{-1}}(h_{(22,p^{-1})}).
    \end{eqnarray*}
    So we conclude that $L_{p} = M_{p}$. This proves both relation \eqref{3.21} and the existence of antipode.
  \end{itemize}
 $\hfill \square$


 Following this theorem, we can get some special results, e.g., Theorem 3.4 in \cite{WL14} and Theorem 3.3 in \cite{JM14}.

 \textbf{Corollary \thesection.4}
 If group-cograded Hopf coquasigroup $H=\bigoplus_{p\in G}H_{p}$ satisfies coassociativity, then the group-cograded Hopf coquasigroup-Ore extension is the Hopf coquasigroup-Ore extension in the sense of \cite{WL14}. The main conclusion Theorem 4.3 in this article is the Theorem 3.4 in \cite{WL14}.
 \\

  \textbf{Corollary \thesection.5}
 If the Group $G$ is trival i.e. $G=\{1\}$ then $R_{1}=H_{1}[y_{1}; \tau_{1}, \delta_{1}]\}$ is the Hopf coquasigroup-Ore extension in the sense of \cite{JM14}. The main conclusion Theorem 4.3 in this article is the Theorem 3.3 in \cite{JM14}.
\\

 \textbf{Proposition \thesection.6}
 If the group-cograded Hopf coquasigroup $R = \bigoplus_{p\in G}R_{p} = \bigoplus_{p\in G}H_{p}[y_{p}; \tau_{p}, \delta_{p}]$ is the group-cograded Hopf coquasigroup-Ore extension, then
 \begin{eqnarray*}
      \Delta_{p,q}(\delta_{pq}(r_{p,q})r_{p,q}^{-1}) = \delta_{p}(r_p)r_p^{-1}\otimes 1 + r_{p}\otimes \delta_{q}(r_q)r_q^{-1}.
 \end{eqnarray*}

 \emph{Proof}
 Since $\Delta_{p,q}(r_{pq}^{1}) = r_{p}^{1}\otimes r_{q}^{1}$ and $r_p = r_{p}^{1}(r_{p}^{2})^{-1}$, it follows that
 \begin{eqnarray*}
      \Delta_{p,q}(r_{p,q}) &=& \Delta_{p,q}(r_{pq}^{1}(r_{pq}^{2})^{-1})\\
      &=& \Delta_{p,q}(r_{pq}^{1})\Delta_{p,q}((r_{pq}^{2})^{-1})\\
      &=& (r_p^1\otimes r_q^1)((r_{p}^{2})^{-1}\otimes (r_{q}^{2})^{-1})\\
      &=& r_{p}^{1}(r_{p}^{2})^{-1}\otimes r_{q}^{1}(r_{q}^{2})^{-1}\\
      &=& r_p\otimes r_q.
\end{eqnarray*}
According to formula \eqref{D3}, it can be obtained
\begin{eqnarray*}
      \Delta_{p,q}(\delta_{pq}(r_{p,q})r_{p,q}^{-1})
      &=& \Delta_{p,q}(\delta_{pq}(r_{p,q}))\Delta_{p,q}(r_{p,q}^{-1})\\
      &=& (\delta_p(r_p)\otimes r_q + r_p^2\otimes \delta_q(r_q))(r_p^{-1}\otimes r_q^{-1})\\
      &=& \delta_{p}(r_p)r_p^{-1}\otimes 1 + r_{p}\otimes \delta_{q}(r_q)r_q^{-1}.
\end{eqnarray*}
$\hfill \square$
\\

 At the end of this section, let $H$ and $H'$ are group-cograded Hopf coquasigroups,
 $R = \bigoplus_{p\in G}R_{p} = \bigoplus_{p\in G}H_{p}[y_{p}; \tau_{p}, \delta_{p}]$ be the Ore extension of $H$ with $\Delta_{p,q}(y_{pq}) = y_{p} \otimes 1_{q} + r_{p} \otimes y_{q}$,
 and $R' = \bigoplus_{p\in G}R'_{p} = \bigoplus_{p\in G}H'_{p}[y'_{p}; \tau'_{p}, \delta'_{p}]$ be the Ore extension of $H'$ with $\Delta_{p,q}(y'_{pq}) = y'_{p} \otimes 1_{q} + r'_{p} \otimes y'_{q}$.
 Now we discuss the isomorphism between Ore extensions of two group-cograded Hopf coquasigroups.
 \\

 \textbf{Theorem \thesection.7}
 Let $R $ and $R' $  be the group-cograded Hopf coquasigroup-Ore extension of $H$ and $H'$ respectively.
 Then $R'$ is isomorphic to $R$ if there is an isomorphism $\varphi: H\rightarrow H'$ such that
 \begin{eqnarray*}
      &&\varphi_p(r_p) = r'_p,\quad \tau'_p(\varphi_p(h_p)) = \varphi_p(\tau_p(h_p)),\\
      &&\delta'_p(\varphi_p(h_p)) = \varphi_p(\delta_p(h_p)) + \varphi_p(\tau_p(h_p))d_p - d_p\varphi_p(h_p),
 \end{eqnarray*}
 where $d_p\in H'_p$ such that $\Delta_{p,q}(d_{pq}) = d_{p} \otimes 1 + r'_{p} \otimes d_{q}$.

 \emph{Proof}
 Let $\overline{\varphi}_p(y_p) = y'_p + d_p$. Then $\varphi$ can be extended from $\varphi: H\rightarrow H'$ to $\overline{\varphi}: R\rightarrow R'$. For all $h_p\in H_p$ we have
 \begin{eqnarray*}
      \overline{\varphi}_p(y_ph_p) &=& \overline{\varphi}_p(\tau_p(h_p)y_p + \delta_p(h_p)) \\
      &=& \overline{\varphi}_p(\tau_p(h_p))\overline{\varphi}_p(y_p) + \overline{\varphi}_p(\delta_p(h_p)) \\
      &=& \overline{\varphi}_p(\tau_p(h_p))(y'_p + d_p) + \overline{\varphi}_p(\delta_p(h_p))
      \\
      \overline{\varphi}_p(y_p)\overline{\varphi}_p(h_p)
      &=& (y'_p + d_p)\overline{\varphi}_p(h_p) \\
      &=& y'_p\overline{\varphi}_p(h_p) + d_p\overline{\varphi}_p(h_p) \\
      &=& \tau'_p(\overline{\varphi}_p(h_p))y'_p + \delta'_{p}(\overline{\varphi}_p(h_p)) + d_p\overline{\varphi}_p(h_p) \\
      &=& \tau'_p(\overline{\varphi}_p(h_p))y'_p + \overline{\varphi}_p(\delta_p(h_p)) + \overline{\varphi}_p(\tau_p(h_p))d_p - d_p\overline{\varphi}_p(h_p)  + d_p\overline{\varphi}_p(h_p) \\
      &=& \tau'_p(\overline{\varphi}_p(h_p))(y'_p + d_p) + \overline{\varphi}_p(\delta_p(h_p)) \\
      &=& \overline{\varphi}_p(\tau_p(h_p))(y'_p + d_p) + \overline{\varphi}_p(\delta_p(h_p)) \\
 \end{eqnarray*}
 So we obtain
 \begin{eqnarray*}
      \overline{\varphi}_p(y_ph_p) = \overline{\varphi}_p(y_p)\overline{\varphi}_p(h_p)
 \end{eqnarray*}
 Since
 \begin{eqnarray*}
      \Delta_{p,q}(\overline{\varphi}_{pq}(y_{pq})) &=&  \Delta_{p,q}(y'_{pq} + d_{pq})\\
      &=& \Delta_{p,q}(y'_{pq}) + \Delta_{p,q}(d_{pq})\\
      &=& y'_{p} \otimes 1_{q} + r'_{p} \otimes y'_{q} + d_{p} \otimes 1_{q} + r'_{p} \otimes d_{q}\\
      &=& (y'_{p} + d_{p}) \otimes 1_{q} + r'_{p} \otimes (y'_{q} + d_{q}),
      \\
      (\overline{\varphi}_p \otimes \overline{\varphi}_q)\Delta_{p,q}(y_{pq}) &=& (\overline{\varphi}_p \otimes \overline{\varphi}_q)(y_{p} \otimes 1_{q} + r_{p} \otimes y_{q})\\
      &=& \overline{\varphi}_q(y_{p}) \otimes 1_{q} + \overline{\varphi}_p(r_{p})\otimes \overline{\varphi}_q(y_{p})\\
      &=& (y'_{p} + d_{p}) \otimes 1_{q} + r'_{p} \otimes (y'_{q} + d_{q}),
\end{eqnarray*}
it follows that
\begin{eqnarray*}
      \Delta_{p,q}(\overline{\varphi}_{pq}(y_{pq})) = (\overline{\varphi}_p \otimes \overline{\varphi}_q)\Delta_{p,q}(y_{pq}).
\end{eqnarray*}
Then we have $R \cong R'$.
$\hfill \square$

\end {document}